\documentclass[12pt]{svproc}
\usepackage{amsmath,amssymb,url}
\usepackage{algorithmic}
\usepackage{graphicx}
\usepackage{enumitem}

\newcommand\subdir{\setbox0=\hbox{$\mathchar "0362$}%
\mathop{\displaystyle\smash{\raise 0.3\ht0\hbox{\hbox to \wd0%
{\hfil\vrule height4pt width0.4pt\hfil}\kern-\wd0%
\lower\ht0\box0}}}
}

\newcommand\Rad{\mbox{\rm Rad}}
\newcommand\Soc{\mbox{\rm Soc}}
\newcommand\Aut{\mbox{\rm Aut}}

\newcommand\gen[1]{\langle#1\rangle}
\newcommand\sz[1]{\left|#1\right|}

\newcommand{\dnk}[3]{\lower 0.2ex\hbox{$#1$}\kern -.1em{\setminus}
                     \kern -.1em\raise 0.2ex\hbox{$#2$}
                     \kern -.1em{/}\kern -.1em\lower 0.2ex\hbox{$#3$}}

\newtheorem{project}{Open Problem}

\begin{document}
\mainmatter
\title{Calculating Subgroups with {\sf GAP}}
\author{Alexander Hulpke}
\institute{Department of Mathematics, Colorado State University, 1874 Campus
Delivery, Fort~Collins, CO, 80523-1874, USA,\\
\email{\url{hulpke@colostate.edu},\\
\url{http://www.math.colostate.edu/~hulpke}}}

\maketitle
\begin{abstract}
We survey group-theoretic  algorithms for finding (some or all) subgroups
of a finite group and discuss the implementation of these algorithms in the computer algebra system {\sf GAP}.
\end{abstract}

One of the earliest questions posed for the development of group theoretic
algorithms has been the determination of the subgroups of a finite group
$G$, as well as the associated lattice
structure.

Since $G$ acts on its subgroups an obvious storage improvement is to store
the subgroups as conjugacy classes, representing each class by a subgroup
$U$ and a transversal of coset representatives of $N_G(U)$ in $G$.

The purpose of this article is to survey the methods that are currently in
use for such computations, not with an aim to supersede the original
descriptions~\cite{neubueser60,hulpkeinvar,coxcannonholtlatt,hulpketfsub} or
to give an implementable description, but to given an
overview of the methods employed. This should allow the reader to understand
the interplay of the methods employed, computational tools required, scope
of calculations, and potential for adaption or modifications by users.

On the way we will indicate a number of open problems, whose solution would
lead to improvements of theoretical or practical aspects of the algorithms.
\medskip

While we shall point to the {\sf GAP} functions that implement the
respective functionality, we shall stop short of printing transcripts of
system sessions, instead the reader is referred to the system documentation.

Neither is this paper intended as a complete survey of Computation Group theory
over its history of at least 60 years. We thus do not aim to cite every
relevant work, but give preference to handbooks or summary articles that are
often easier accessible.
\smallskip

We will illustrate the scope of calculations by assuming a contemporary
(as of 2017) standard desktop machine with a 3.5GHz processor (utilizing just a
single core) and 8GB of memory.

\section{Tools Required}

In general we will represent a subgroup $S$ of the finite group $G$ by a set
of generators, given as elements of $G$. One may think of $G$ as the group containing all
transformations of a given kind --- for example in the case of permutations
a symmetric group $S_n$ or even the finitary symmetric group on positive
integers. Similarly, in the case of matrices this group might be the full
general linear group.

We thus need methods that allow us to determine for such a
subgroup $S$:
\begin{itemize}
\item The order of $S$.
\item Test whether an element of $G$ is contained in $S$, and if so:
\item Express an element of $S$ as a word in the given generators of $S$,
thus enabling us to
evaluate homomorphisms.
\item Write a presentation for $S$ in a given generating set, thus
testing whether a map on generators is a homomorphism. (In practice one often
does not use an arbitrary generating set, but a specific one that allows for
a nicer presentation.)
\item Determine a composition series, a chief series and the radical
$\Rad(G)$ (the largest solvable normal subgroup) of $G$, as well as a
representation of $G/\Rad(G)$ as a permutation or matrix group.
\end{itemize}

For permutation groups, such functionality is obtained through a stabilizer
chain~\cite[Chapter 4]{holtbook}, respectively~\cite{akosbook}.
For matrix groups, such functionality is
provided by the data structure of a composition
tree~\cite{baarnhielmholtleedhamobrien,neunhoefferseress06}.
These tools can be extended to groups of other classes of invertible
transformations of a finite object using the ``black-box''
paradigm~\cite{babaibealsseress09}.

For solvable groups, polycyclic generating sets
(that is a set
of generators that is adapted to a composition series and allows for an
effective normal form) provide such functionality~\cite{SOGOS}, see
also~\cite[Chapter 8]{holtbook}.

\subsection{Complexity}
For solvable groups, polynomial time algorithms are known for all of these
tasks.

For permutation groups, the known algorithms are proven to be
polynomial time, as long as no composition factor of type ${}^2G_2(q)$
occurs (in which case the result will still be correct, but the time bound is not known to
hold.). In fact, the algorithms are almost linear (linear up to logarithmic
factors) time in a Las Vegas probabilistic setting (see \ref{randel} below).

\begin{project}
Show that the groups ${}^2G_2(q)$ have a short presentation in the sense
of~\cite{babaigoodmanetal97}. Such a result will allow the removal of the
qualifier in the previous paragraph.
\end{project}

The complexity situation for matrix groups~\cite{babaibealsseress09} is as
with permutation groups with one further complication: $\mbox{GL}_n(q)$ contains
cyclic subgroups (Singer cycles) of order $q^n-1$, and calculations in these
groups are equivalent to Discrete Logarithm problems. The proven complexity
is therefore also up to a Discrete Logarithm "oracle", that is the cost of
discrete logarithm calculations is not accounted for.

These polynomial time algorithms for solvable and for permutation groups
have been fully implemented in {\sf GAP} and in {\sf Magma}. The available
implementations for the matrix group algorithms involve many, but not all,
of the polynomial time methods.  The reason for this is that there are
number of algorithms for subtasks that perform better in practice than the
generic black-box algorithm, but so far no proof of polynomial time has been
found.

\medskip

Arbitrary finitely presented groups will require the use of a faithful
representation in one the respresentations discussed before.

\subsection{Random Elements}
\label{randel}
Some of the algorithms utilize random selections of elements. It thus seems
appropriate to briefly address this issue.

First, on the computer random selection is always based on a random number
generator, and thus is inherently pseudo-random.

Secondly, once we can test membership of elements, the underlying data
structures allow us to construct a bijection between $G$ and the numbers
$1,\ldots,{\sz{G}}$ and thus select elements of the same random quality as
the random number generator provides.
\medskip

Some of the functions to build basic data structures also utilize 
pseudo-random elements which are obtained as pseudo-random products of
generators and inverses~\cite{celleretalrandom,babaipak04}.
All of these calculations then involve verification steps that ensure the
returned result is always correct, regardless of the random choices or the
quality of randomness.

As far as complexity is concerned, any such algorithm then lies in a class
denoted by ``Las Vegas'': That is the algorithm will always return a correct
result and will, with a user-chosen
probability $0<\varepsilon<1$, terminate in the given time.
However with probability $1-\varepsilon$ the calculation will take longer (but will
eventually terminate with a correct result).

\subsection{Mid-level tools}
Building on these tools, a number of mid-level tools obtain
structural group-theoretic information:
\begin{itemize}
\item For $S\le G$, representatives of the cosets of $S$ in
$G$ \cite{dixonmajeed}.
\item The centralizer $C_G(g)$ of elements $g\in G$ as well as conjugating
elements $x$ that for given $g,h\in G$ satisfy $g^x=h$ (if they exist). (For
permutation groups this is a backtrack search, following~\cite{leon91}).
\item The normalizer $N_G(S)$ of a subgroup $S\le G$ as well as conjugating
elements $x$ that for given $S,T\le G$ satisfy that $S^x=T$ \cite{leon97}.
\item Representatives of the conjugacy classes of elements of
$G$
\cite{meckyneubueser,cannonsouvignier,hulpkepermclass,cannonholt07,hulpkematclass}.
\item Representatives of Sylow subgroups of $G$ for a chosen prime.
\item For a normal subgroup $N\lhd G$, representatives of the $G$-classes
of complements to $N$ in $G$, provided that $N$ is
solvable~\cite{cellerneubueserwright}.
This algorithm is based on cohomomology through a presentation for $G/N$.

If $G/N$ is solvable, complements can be computed in a combination
of cohomology and reduction to subgroups \cite{hulpketfsub}.
\item
Determine an effective\footnote{Meaning that it, and its inverse can be
applied to group elements to obtain the image}
isomorphism between two groups $G$ and $H$ (or show
that no such isomorphism can exist)~\cite{obrien93,holtcannonautgroup}.
\end{itemize}
These algorithms are typically not of polynomial, but exponential worst case
time complexity. However in most cases of practical interest they tend to work
well, 
allowing for them to be used as building blocks for larger calculations.

\section{The Basic Structure}
\label{bastru}

\begin{figure}
\begin{center}
\includegraphics[width=7cm]{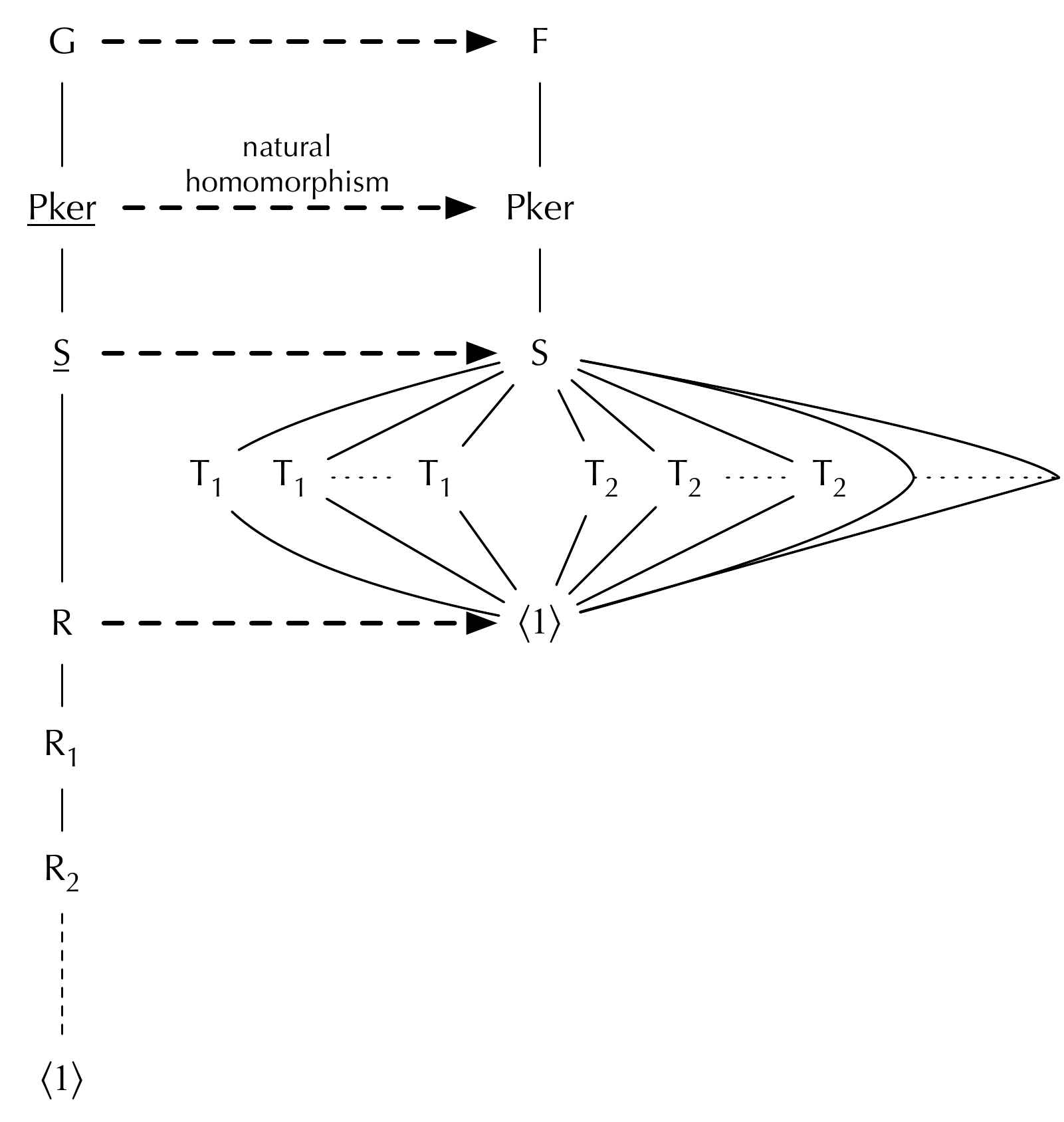}
\end{center}
\caption{Subgroups related to the solvable radical data structure}
\label{solvrad}
\end{figure}

The basic structure underlying most subgroup calculations and the one we shall use is based on the {\em solvable radical} (or
{\em trivial fitting}) paradigm~\cite{babaibeals99,holt97,cannonsouvignier}, as
depicted in Figure~\ref{solvrad}:

Let $G$ be a finite group, $R=\Rad(G)$ and $\varphi\colon G\to G/R=:F$. Then
$S=\Soc(F)=\prod T_i$ must be the direct product of nonabelian simple groups
$T_i$. We thus
can assume that $F$ is represented as a subgroup of $\Aut(\Soc(F))$; that is
as a subgroup of a direct product of groups of the form $\Aut(T_i)\wr
S_{m_i}$ for
$T_i$ simple and $\sum_i m_i $ the number of simple factors of $S$.

The action of $F$ on the socle factors has a kernel denoted by $Pker$, the
factor $Pker/S$ is a direct product of subgroups of outer
automorphisms. We denote by $\underline{S}$ and $\underline{Pker}$ the full
preimages of these subgroups in $G$.

We now determine subgroups in the following way.
\begin{enumerate}
\item
\label{oans}
Subgroups of the simple socle factors $T_i$.
\item
\label{zwoa}
Combine these to subgroups of $\Soc(F)$.
\item
Calculate the subgroups of $F/\Soc(F)$ (which will be a much smaller group
than $F$).
\item
\label{gsuffa}
Extend the subgroups of $\Soc(F)$ to subgroups of $F$ by using the subgroups
of $F/\Soc(F)$.
\item Determine a series of normal subgroups
$R=R_0>R_1>R_2>\cdots>R_k=\gen{1}$ with $R_i\lhd G$ and $R_i/R_{i+1}$
elementary abelian.
\item
\label{gsuffa2}
Determine subgroups of $G/R_{i+1}$ from subgroups of $G/R_i$
(initialized for $i=0$ with $G/R_0=F$) and the $G$-module action on
$R_i/R_{i+1}$. Iterate.
\end{enumerate}

Typically we will store not all subgroups of a group $G$, but only
representatives of the conjugacy classes under $G$, since this saves
substantially on the memory requirements.
This enumeration up to conjugacy can be translated
for each of these steps to conjugacy under suitable actions. For example
in step~\ref{oans} it is conjugacy by the subgroup of $\Aut(T_i)$ induced
through the action of $N_F(T_i)$. Finding representatives up to conjugacy
can in general mean that we have to do explicit subgroup conjugacy tests.
In some steps of the algorithm  (say when calculating complements by
cohomological methods) such tests can be preempted or reduced by using 
other equivalences amongst the objects constructed.

In the following more detailed description we shall focus on the
task of finding all groups rather than the elimination of conjugates.
\medskip

Methods similar to section~\ref{include} can then be used to determine the
incidence structure of the full subgroup lattice.

\section{The steps of the algorithm}

We now describe the different steps of the algorithm in more detail:

\subsection{Factor Groups}

A fundamental paradigm of the approach is to work in homomorphic images.
This raises the question of how to represent factor groups of $G$ in a
suitable way. While this is difficult in general, for the particular factor
groups required here effective solutions exist:

\begin{itemize}
\item
It has been shown~\cite{LuksSeress97,holt97} that for permutation groups
$G$, the factor $G/\Rad(G)$ can be (constructively) represented with
permutation degree not exceeding that of $G$. (In {\sf GAP} this is a call
to {\tt NaturalHomomorphismByNormalSubgroup(G,RadicalGroup(G))}.

More generically, the special structure of $G/\Rad(G)$ as a subgroup of a
direct product of wreath products allows for a representation of moderate
degree, using imprimitive wreath products.
\item By Schreier's conjecture (as proven in~\cite{Feit80}), for a simple group
$T$ the factor $\Aut(T)/T$ is small. Thus $F/\Soc(F)$ (which embeds
into a direct product of groups of the form $(\Aut(T_i)/T_i)\wr S_{m_i}$) is
comparatively small and can be easily represented in an ad-hoc way.
\item
In many cases it is not necessary to represent a factor group $G/N$
faithfully, but it is sufficient to use representatives of elements and full
preimages of subgroups. In particular, we can use this to perform linear
algebra with coefficient vectors for the abelian factors
$R_i/R_{i+1}$ of the radical.
\end{itemize}

The question of the minimal permutation degree of factor groups of
permutation groups has been studied also theoretically, and one can ask for
other classes of normal subgroups for which such degree bounds hold:
\begin{project}
Extending the work of~\cite{easdownpraeger}, describe (constructively) cases
in which for permutation groups or matrix groups $G$ and $N\lhd G$ one can represent
the factor group $G/N$ in degree not exceeding that of $G$. 
\end{project}

\subsection{Subgroups of simple groups}

Step~\ref{oans} (from page~\pageref{oans}) asks us to determine the subgroups of a simple group $T$.

The basic method for this is the ``cyclic extension'' algorithm, dating back
to~\cite{neubueser60}:
A subgroup $S\le T$ is either perfect, or there is a smaller subgroup $S'\le
U<S$ such that $S=\gen{U,n}$ with $n\in N_G(U)$.
Thus:
\begin{itemize}
\item[a)] Initialize the perfect subgroups of $T$. This requires a
precomputed list of isomorphism types of perfect groups such
as~\cite{holtplesken91} for groups of order at most $10^6$.
(By now, due to the rapid progress in computer
engineering, the same methods would allow us to build such lists for larger
orders.)

Then, in an approach close to isomorphism test algorithms, search for
isomorphic copies of each of these groups as subgroups of $T$.

In {\sf GAP} such a list is obtained using the operation {\tt
RepresentativesPerfectSubgroups}.
\item[b)] For every subgroup $U$ listed so far, classify the $U$-orbits of
elements of $N_G(U)$ outside $U$. If for an orbit representative $n$ the group
$\gen{U,n}$ is not yet known (i.e. not conjugate to a known group) then add it to the list. Iterate.
\end{itemize}
To allow for an efficient storage/comparison of subgroups, the algorithm
maintains a list of cyclic subgroups of prime power order (called zuppos by
their German acronym\footnote{``Zyklische Untergruppen von Primzahl-Potenz
Ordnung''}). It then represents every subgroup as a bit list indicating
which zuppos it contains.
\medskip

Simple groups tend to have relatively few subgroups, enabling the
calculation of subgroups even for large group orders. The assumed standard
computer will calculate the subgroups of a simple group of order $10^5$ in
under a minute, order $10^6$ about 5-10 minutes and (provided the potential
perfect subgroups are available) order $10^7$ about 90 minutes. (This is
assuming that the group is given as a permutation group of minimal degree.)

The algorithm of course also will work for groups that are not simple, but
in this case is often not competitive.
\medskip

In {\sf GAP}, this algorithm is implemented by the command {\tt
LatticeByCyclicExtension}.

\bigskip

In practice, we can (using this algorithm) create a database
of subgroups of simple groups
$T$ up to a certain order limit once, and then store them.
If the algorithm then is called for one of these
simple groups, one
then simply can fetch subgroups from the database.
\smallskip

{\sf GAP} does exactly this, the databases used to obtain subgroup
information is the library of tables of marks, provided by the {\tt tomlib}
package (which will be loaded automatically, if available). 
As of writing, this library contains full subgroup data for most of the
simple groups in the ATLAS of order roughly up to $10^7$.
Some information
about maximal subgroups of symmetric and alternating groups is also obtained 
through the library of primitive groups.

\bigskip

This approach requires an isomorphism between the concrete simple group $T$ and
its incarnation $D$ in the database. Such an isomorphism can be facilitated
in many cases through the use of so-called {\em standard
generators}~\cite{wilson98}: For a simple group $T$, this is a pair of
elements $a,b\in T$ such that:
\begin{itemize}
\item $T=\gen{a,b}$. (By~\cite{aschbacherguralnick84} every finite simple
group can be generated by two elements.)
\item
The pair $(a,b)$ (that is its $\Aut(T)$-orbit) is characterized
by simple relations, such as orders of $a$ and $b$ or short product
expressions in $a$ and $b$, or $T$-class memberships of $a$ and $b$.
This implies
that if $T_1\cong T_2\cong T$ an isomorphism $T_1\to T_2$ is obtained by
finding instances of standard generators $a_1,b_1\in T_1$ and $a_2,b_2\in
T_2$ and constructing
the homomorphism that maps $a_1$ to $a_2$ and $b_1$ to $b_2$.
\item
In a given instance of $T$, such a pair $(a,b)$ can be found quickly by only
using basic group operations such as product and inverse (thus allowing for
pseudo-random elements) and element order. A typical property achieving
this is if the elements lie in small conjugacy classes that are powers of
large conjugacy classes: A (pseudo-)random element will likely lie in a large class,
by powering we get an element in the small class and only few conjugates to
consider.

\end{itemize}
For example $\sz{a}=2$, $\sz{b}=3$, $\sz{ab}=5$ could be used as such a
generating set for $A_5$.

Such standard generators have been defined for all sporadic groups and many
groups of Lie type of small order.
\begin{project}
Generalize ``standard generators'' to all quasisimple groups of Lie type.
\end{project}
\medskip

The concept of standard generators can be generalized to
{\em constructive recognition}, that is the task to find an isomorphism from
a simple group $T$ to its stored database incarnation $D$, without relying
on the need to find specific generators, but rather ``rebuilding'' natural
combinatorial structures from within the group. For example, if the group
$T$ is a matrix group isomorphic to $A_n$, one might want to find a subspace
of the natural module that has an orbit of length $n$ under $T$, thus
providing such an isomorphism through the action on the subspaces in the
orbit.
See the survey~\cite{dietrichleedhamgreenobrien} for formal
definitions and details.

\subsection{Subdirect products}

Step~\ref{zwoa} combines the subgroups of direct factors to those of a
direct product. By induction it is sufficient to consider the case of a
direct product of two groups, $G\times H$. Let $S\le G\times H$ and denote
the projection from $S$ to $G$ by $\alpha$ and that from $S$ to $H$ by
$\beta$. The image groups $A=S^\alpha$ and $B=S^\beta$ then are subgroups of
$G$, respectively $H$.

Given such subgroups $A$ and $B$, the construction of a {\em subdirect
product} (which dates back at least to~\cite{remak})
then allows to construct all groups $S$ (see Figure~\ref{subdir}):

\begin{figure}
\begin{center}
\includegraphics[width=10cm]{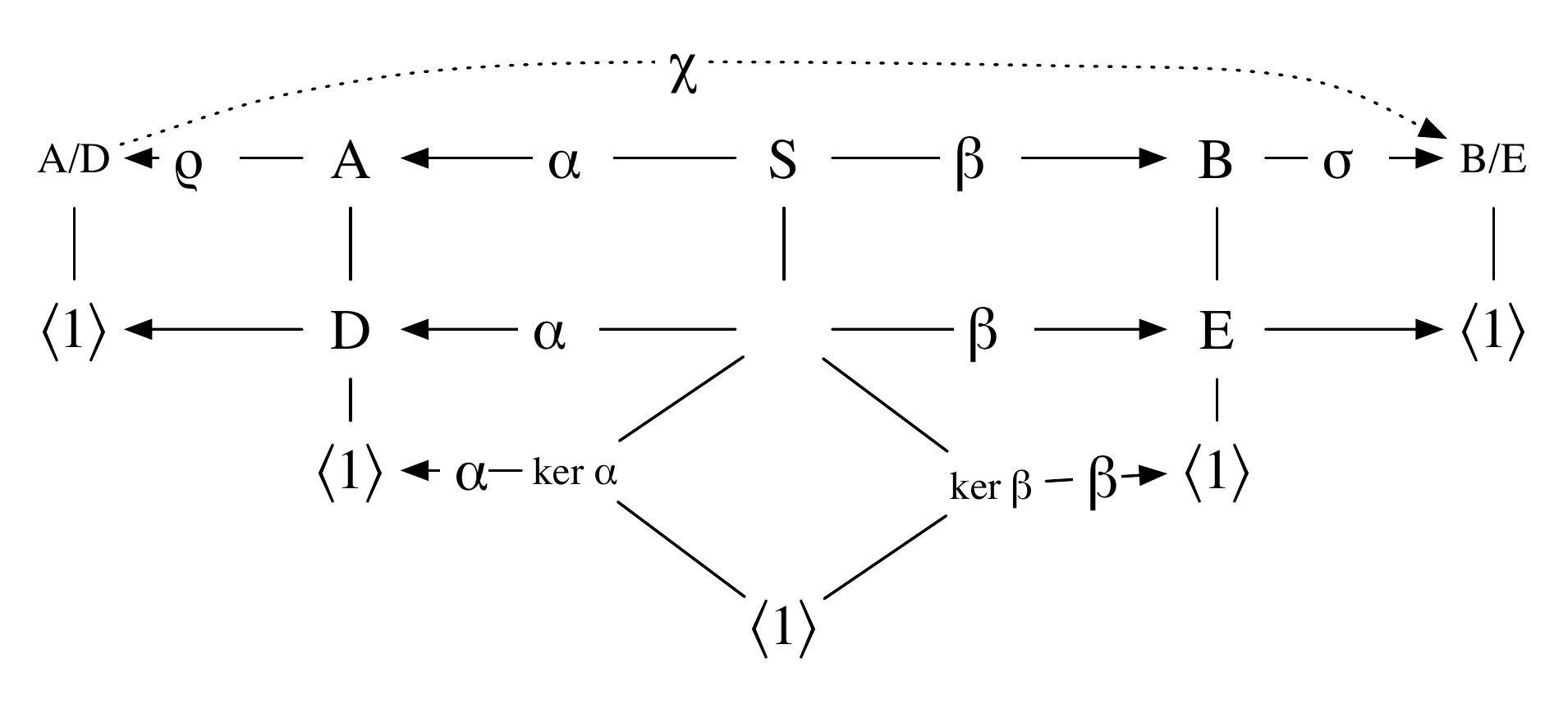}
\end{center}
\caption{Subdirect product construction}
\label{subdir}
\end{figure}

Denote by $D\lhd A$ the image of $\ker\beta$ under $\alpha$ and by $E\lhd B$
the image of $\ker\alpha$ under $\beta$. Then by the isomorphism theorem
\[
A/D\cong S/\gen{\ker\alpha,\ker\beta}\cong B/E.
\]
If $\chi\colon A/D\to B/E$ is this isomorphism, and we denote the natural
homomorphisms by $\varrho\colon A\to A/D$ and $\sigma\colon B\to B/E$, then
\[
S=\left\{(a,b)\in G\times H\mid a\in A, b\in B,
(a^\varrho)^\chi=b^\sigma\right\}.
\]

To construct all subdirect products $S$ corresponding to the pair $A,B$,
we thus classify pairs of normal subgroups $D\lhd A$, $E\lhd F$ together
with isomorphisms $\chi\colon A/D\to B/E$.

Conjugacy of subgroups by 
$N_G(A)\times N_G(B)$ will induce equivalences on the normal subgroups and
amongst the isomorphisms.
\smallskip

In the case we consider -- subgroups of $\Soc(F)$ -- furthermore there may
be a conjugation action of $F$ on the direct factors of its socle that
causes further fusion of subgroups.

\subsection{Normal Subgroups and Complements}

In steps~\ref{gsuffa} and~\ref{gsuffa2} of the calculation, we have a normal
subgroup $N\lhd G$ and know the subgroups of $G/N$ as well as the subgroups of
$N$. (In step~\ref{gsuffa2} the normal subgroup $N$ is a
vector space whose subgroups are easily enumerated.)
From these we want to construct the subgroups of $G$.

\begin{figure}
\begin{center}
\includegraphics[width=7cm]{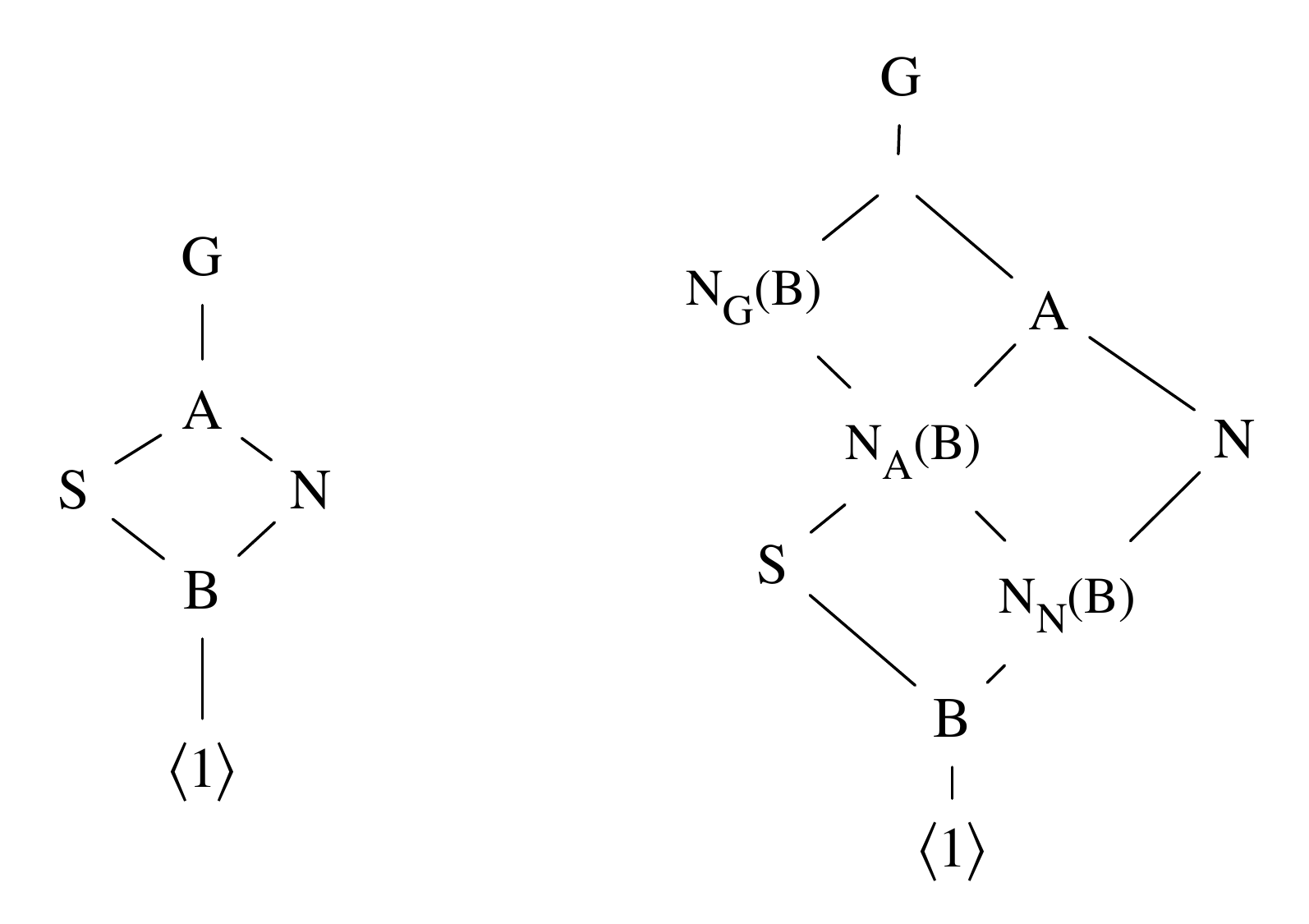}
\end{center}
\caption{Complement situations for subgroups}
\label{compgrp}
\end{figure}

We first analyze the situation: Let $S\le G$ and set $A=\gen{N,S}$ and
$B=S\cap N\lhd B$.
(See Figure~\ref{compgrp}, left.)

\paragraph{A) Abelian Normal subgroup}
We consider first the case that $N$ is abelian (which arises in
step~\ref{gsuffa2}). Then $B\lhd N$ and thus $B\lhd\gen{S,N}=A$.

Thus $S/B$ is a complement to $N/B$ in $A/B$. As
$N/B$ is elementary abelian, such complements can be obtained through
cohomological methods, following~\cite{cellerneubueserwright}.
The input to such a computation is the linear action of $A$ on
$N/B$, together with a presentation for $A/N$.
\smallskip

To find all subgroups of $G$, we iterate through all $A$ (as subgroups of
$G/N$) and for each $A$ determine candidates for $B$ as submodules of $N$
under the action of $A$~\cite{luxmuellerringe}.

\smallskip
As step~\ref{gsuffa2} then iterates over a series, a crucial step towards
efficiency is to extend a presentation for $A/N$ to a presentation for $S$,
if $S/B$ is such a complement. This is easy, as $B$ is elementary
abelian.

\paragraph{B) Nonabelian Normal subgroup}
If $N$ is not abelian (as it will be in step~\ref{gsuffa}), the situation is
more complicated, as $B$ is not necessarily normal in $A$, and there is no
algorithm to easily determine complementing subgroups. In this
case, following~\cite{hulpketfsub}, we
iterate through the possible subgroups $B\le N$ and for each such $B$
determine the groups $S$ such that $S\cap N=B$:

As $N\le\gen{N,S}=A$, we have that $N_N(B)\le N_A(B)=\gen{S,N_N(B)}\le
N_G(B)$. Furthermore, $N_G(B)/N_N(B)$ is isomorphic to a subgroup of $G/N$.
(See Figure~\ref{compgrp}, right.)
In this situation $S/B$ is a complement to $N_N(B)/B$ in $N_A(B)/B$.

Given a subgroup $B\le N$, we thus determine the subgroups of $N_G(B)/N_N(B)$
(e.g. from the subgroups of $G/N$) and for each subgroup $N_A(B)/N_N(B)$
determine the candidates for $S/B$ as complements. If $N_N(B)/B$ is
solvable, this again can be done using cohomology calculations.
\medskip

The group $N_N(B)/B$ does not need to be solvable -- if the factor
group however is solvable (which will be the case unless $\Soc(F)$ contains
a single simple factor at least quintuply, in which case there will be
storage problems already for the subgroups of $\Soc(F)$),
\cite{hulpketfsub} describes an approach for complements that reduces to
$p$-groups, corresponding to a chief series of the factor.
\medskip

{\sf GAP} contains a function {\tt ComplementClassesRepresentatives({\em
G},{\em N})} that determines representatives of the classes of complements
to $N$ in $G$, up to conjugacy by $G$, provided that $N$ or $G/N$ are solvable.

In the case that neither $G$ and $G/N$ are solvable, no algorithm for
complements exists yet:

\begin{project}
Find a good algorithm for determining complements if both normal subgroup
and factor groups are not solvable. This also
has relevance to maximal subgroup computations~\cite{holtcannonmaxgroup}.
\end{project}

\subsection{Implementation}
In {\sf GAP} the algorithm described in the previous sections (with some variants depending on the
representation of the groups) is obtained through the operation {\tt
ConjugacyClassesSubgroups}. It takes as argument a group and returns a list
of conjugacy classes of subgroups. For each class {\tt Representative} will
return one subgroup; {\tt AsList} applied to a class will return all
subgroups in this class, thus
\begin{verbatim}
Concatenation(List(ConjugacyClassesSubgroups(G),AsList));
\end{verbatim}
returns all subgroups of a group $G$. In general such an enumeration of all
subgroups is not recommended as it is very costly in terms of memory.
\medskip

It is also possible to visualize the full lattice of subgroups of a group $G$.
For this, the command
\begin{verbatim}
DotFileLatticeSubgroups(LatticeSubgroups(G),"filename.dot");
\end{verbatim}
produces a text file, called {\tt filename.dot} (or whatever file name is
given) that describes the incidence structure of the subgroup lattice in the
{\tt graphviz} format (see~\url{www.graphviz.org} for a description and for
viewer programs for this format. There also are programs to convert this
format into others, e.g. {\tt dot2tex} converts to {\tt TikZ} or {\tt
PSTricks} format.

\begin{figure}
\begin{center}
\includegraphics[width=11cm]{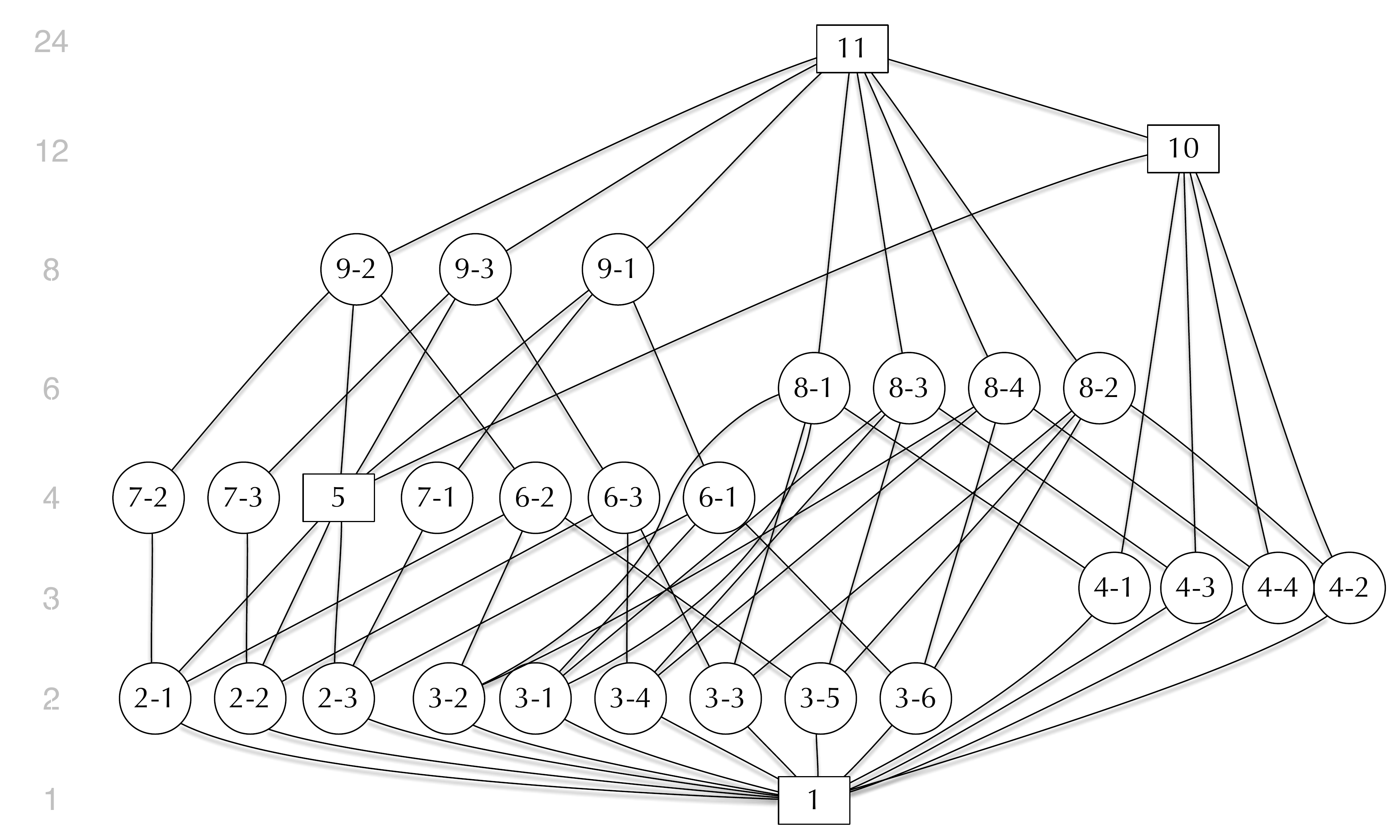}
\end{center}
\caption{Subgroup Lattice for $S_4$}
(produced using {\tt DotFileLatticeSubgroups} and then visualized using the
graphics software {\sf OmniGraffle})
\label{s4}
\end{figure}

Figure~\ref{s4} illustrates the result in
the example of the symmetric group $S_4$. Rectangles represent normal
subgroups, circles ordinary subgroups and their conjugates. A number $a-b$
indicates group number $b$ in class $a$ (there is no $b$-part if the group
is normal, as it will default to $b=1$).

This group can be obtained in {\sf GAP} then as {\tt cl[{\it a}][{\it b}]}
(that is {\tt cl[8][3]} for $a=8$ and $b=3$) where {\tt
cl:=ConjugacyClassesSubgroups(G)}.

{\bf Caveat:} The ordering (both $a$ and $b$-parts) of subgroups can involve
ad-hoc choices within the algorithm. When creating the group $G$ a second
time with the same generators, it is possible that a different numbering is
chosen. It thus is not safe to use the $a-b$ indices for specifying a
concrete subgroup outside a particular run of {\sf GAP}.

\subsection{Practicality and Modifications}

With the construction process proceeding through layers, in each step
proceeding through all subgroups found in the previous step, the limiting
factor to calculation is (as timings in~\cite{hulpketfsub} indicate)
the total number of subgroups, rather than the group order.
\medskip

If only some subgroups are desired, and calculation of the full lattice is
infeasible, it might be possible to restrict the calculations to certain
subgroups, as long as a filter can be defined that is appropriate to the
construction process and will iterate the construction only for subgroups
with certain properties. (For example, the cyclic extension algorithm might
be instructed to not calculate subgroups larger than a prescribed limit.)

At the moment, {\sf GAP} provides
options to define such filters in a few cases (see the manual for details):
\begin{itemize}
\item The general algorithm, as described, is implemented by the operation
{\tt LatticeViaRadical}. If given two groups as argument it calculates
subgroups of the second group up to conjugacy by the first group.
\item {\tt LatticeByCyclicExtension} allows for limiting the extension step
to subgroups with a particular property.
\item {\tt SubgroupsSolvableGroup}, an implementation of the algorithm
described for the case of solvable groups (in which case only
step~\ref{gsuffa2} is needed) allows to limit the determination of
complements to specified cases, depending on properties of $A$, $N$ and $B$.
\item In a different restriction, {\tt SubgroupsSolvableGroup} also allows
for determination of only those subgroups that are fixed (as subgroups)
under a prescribed set of automorphisms, generalizing the concept of
submodules~\cite{hulpkeinvar}.
\end{itemize}

\section{Maximal, Low-Index and Intermediate Subgroups}

A different class of algorithms is obtained by considering maximal
subgroups.

If $M\le G$ is a maximal subgroup, the action of $G$ on the cosets of $M$
is primitive. The classification of primitive groups under the label 
O'Nan-Scott theorem~\cite{onanscott} (see \cite{liebeckpraegersaxl88} for a
full proof with corrections) thus can be used to describe possible maximal
subgroups -- one needs to search for quotient groups of $G$ that have the
correct structure to allow a primitive action, the point stabilizers for
these actions will be maximal subgroups.

An approach to determine representatives
of the conjugacy classes of maximal subgroups of a finite group, using this
idea,
is described in~\cite{eickhulpke01,cannonholt07}. The fundamental
ingredients of these calculations again are the simple factors of $\Soc(F)$,
and complements.
\medskip

Taking again a series as described in Section~\ref{bastru}, the algorithm
then identifies factor groups of $G$ that can have a faithful primitive
representation. This is done via the socle of these subgroups, that is chief
factors (or combinations of chief factors) of $G$:
\begin{itemize}
\item Maximal subgroups intersecting the radical lead to primitive actions of
affine type and thus are obtained as complements. This is the only case of a
solvable socle.
\item 
Nonsolvable chief factors are obtained as part of $\Soc(F)$. Isomorphisms
between the simple factors can be used to construct the different types of
primitive actions, according to the diagonal and product action cases of the
O'Nan-Scott theorem.
\item
The base case is maximal subgroups of simple groups, for which
classifications exist in~\cite{kleidmanliebeck} and (far more
explicitly)~\cite{brayholtdougal}.
\end{itemize}

\begin{project}
Extend the concrete classification of maximal subgroups
in~\cite{brayholtdougal} to larger degrees.
\end{project}

As in the case of using stored tabulated information about subgroups, an
explicit isomorphism needs to be found using constructive recognition or
standard generators.
\medskip

In {\sf GAP}, representatives of the classes of maximal subgroups can be
obtained using the function {\tt MaximalSubgroupClassReps}. (Be aware that
while {\tt MaximalSubgroups} also exists, it will enumerate {\em all}
maximal subgroups, often at significant cost.) Again tabulated information
about maximal subgroups of simple groups is used.

\subsection{Small index and intermediate subgroups}
\label{include}

The maximal subgroup functionality can be used to determine the maximal
subgroups of a subgroup, thus obtaining maximal inclusion.
(This also is used in general to provide the maximality relations required
for the subgroup lattice structure.)

Iterating maximal subgroups can be used to find subgroups that have bounded
index~\cite{cannonholtslatterysteel}, or simply to iterate the computation
of maximal subgroups for all subgroups obtained so far
to find subgroups that are $k$-step maximal in $G$. To reduce the cost it
will be natural to fuse conjugates under the action of the whole group.

In {\sf GAP}, such latter functionality will be provided (starting with the
4.9 release) by a function {\tt LowLayerSubgroups} that for a given group
$G$ and step limit $k$ determines the subgroups of $G$, up to conjugacy,
that are at most $k$-step maximal in $G$. It is possible to limit the
calculation to obtain only subgroups of specified bounded index.
\bigskip

A further variant is to determine the {\em intermediate subgroups} $U<V<G$ for
a given subgroup $U\le G$~\cite{hulpkeintermediate}: Instead of choosing an
arbitrary representative $M$ for each class of maximal subgroups, we
determine in each step which conjugates of $M$ contain the chosen subgroup $U$ and then
iterate.

This variant is implemented in {\sf GAP} by the function {\tt
IntermediateSubgroups} (again this will see a significant performance
improvement with the 4.9 release).

\section{Summary}

We have described the various methods that can be used in {\sf GAP} to
determine the subgroups of a given finite group. Different approaches
provide different options to adapt the calculation. The methods also rely on
a significant framework for basic operations that is essentially invisible to
a user who does not look into the inner workings.
While a calculation of subgroups is mostly limited by the size
of the output set, there are still open research problems whose solution
would improve this (and other) group theoretic algorithms.

\subsection*{Acknowledgment}
The author's work has been supported in part by
Simons Foundation Collaboration Grant~244502.

\bibliographystyle{alpha}

\newcommand{\etalchar}[1]{$^{#1}$}

\end{document}